# STEIN'S METHOD FOR DISCRETE GIBBS MEASURES[1]

By Peter Eichelsbacher and Gesine Reinert

*Ruhr-Universität Bochum and University of Oxford*

Stein's method provides a way of bounding the distance of a probability distribution to a target distribution $\mu$. Here we develop Stein's method for the class of *discrete Gibbs measures* with a density $e^V$, where $V$ is the energy function. Using size bias couplings, we treat an example of Gibbs convergence for strongly correlated random variables due to Chayes and Klein [*Helv. Phys. Acta* **67** (1994) 30–42]. We obtain estimates of the approximation to a grand-canonical Gibbs ensemble. As side results, we slightly improve on the Barbour, Holst and Janson [*Poisson Approximation* (1992)] bounds for Poisson approximation to the sum of independent indicators, and in the case of the geometric distribution we derive better nonuniform Stein bounds than Brown and Xia [*Ann. Probab.* **29** (2001) 1373–1403].

**0. Introduction.** Stein [17] introduced an elegant method for proving convergence of random variables toward a standard normal variable. Barbour [2, 3] and Götze [10] developed a dynamical point of view of Stein's method using time-reversible Markov processes. If $\mu$ is the stationary distribution of a homogeneous Markov process with generator $\mathcal{A}$, then $X \sim \mu$ if and only if $E\mathcal{A}g(X) = 0$ for all functions $g$ in the domain $\mathcal{D}(\mathcal{A})$ of the operator $\mathcal{A}$. For any random variable $W$ and for any suitable function $f$, to assess the distance $|Ef(W) - \int f \, d\mu|$ we first find a solution $g$ of the equation

$$\mathcal{A}g(x) = f(x) - \int f \, d\mu.$$

If $g$ is in the domain $\mathcal{D}(\mathcal{A})$ of $\mathcal{A}$, then we obtain

$$(0.1) \qquad \left| Ef(W) - \int f \, d\mu \right| = |E\mathcal{A}g(W)|.$$

Received April 2006; revised October 2007.

[1]Supported in part by the DAAD Foundation and the British Council, ARC Project, Contract 313-ARC-00/27749 and by the RiP program at Oberwolfach, Germany.

Bounding the right-hand side of (0.1) for a sufficiently large class of functions $f$ leads to bounds on distances between the distributions. Here we will mainly focus on the *total variation distance*, where indicator functions are the test functions to consider.

Such distributional bounds are not only useful for limits, but they are also of particular interest, for example, when the distance to the target distribution is not negligible. Such relatively large distances occur, for example, when only few observations are available, or when there is considerable dependence in the data slowing down the convergence. The bound on the distance can then be taken into account explicitly, for example when deriving confidence intervals.

In Section 1 we introduce Stein's method using the generator approach for *discrete Gibbs measures*, which are probability measures on $\mathbb{N}_0$ with probability weights $\mu(k)$ proportional to $\exp(V(k))$ for some function $V : \mathbb{N}_0 \to \mathbb{R}$. The discrete Gibbs measures include the classical distributions Poisson, binomial, geometric, negative binomial, hypergeometric and the discrete uniform, to name but a few. One can construct simple birth–death processes, which are time-reversible, and which have a discrete Gibbs measure as its equilibrium measure. In the context of spatial Gibbs measures this connection was introduced by Preston [15].

In Section 2 we not only recall bounds for the increments of the solution of the Stein equation from [5], but we also derive bounds on the solution itself, in terms of potential function of the Gibbs measure; see Lemmas 2.1, 2.4 and 2.5. The bounds, which to our knowledge are new, are illustrated for the Poisson, the binomial and the geometric distribution.

For nonnegative random variables, the size bias coupling is a very useful approach to disentangle dependence. Its formulation and its application to assess the distance to Gibbs measures are described in Section 3. We compare the distributions by comparing their respective generators, an idea also used in [11], while paying special attention to the case that the domains of the two generators are not identical. The size bias coupling then naturally leads to Theorem 3.5 and Corollary 3.8.

Section 4 applies these theoretical results to assess the distance to a Gibbs distribution for the law of a sum of possibly strongly correlated random variables. The main results are Theorem 4.2 and Proposition 4.4, where we give general bounds for the total variation distance between the distribution of certain sums of strongly correlated random variables and discrete Gibbs distributions of a grand-canonical form; see [16], Section 1.2.3. In particular, Theorem 4.2 gives a bound on the rate of convergence for the qualitative results in [6] by bounding the rate of convergence.

Considering two examples with nontrivial interaction, we obtain bounds to limiting nonclassical Gibbs distributions. Our bound on the approximation error is phrased in terms of the particle number and the average density of the particles.

Summarizing, the main advantage of our considerations is the application of Stein's method to models with interaction described by Gibbs measures. When applying our bounds on the solution of the Stein equation to the Poisson distribution and the geometric distribution, surprisingly we obtain improved bounds for these well-studied distributions. Thus our investigation of discrete Gibbs measures serves also as a vehicle for obtaining results for classical discrete distributions.

The results presented here will also provide a foundation for introducing Stein's method for spatial Gibbs measures and Gibbs point processes, in forthcoming work.

## 1. Gibbs measures and birth–death processes.

### 1.1. *Birth–death processes.*

A birth–death process $\{X(t), t \in \mathbb{R}\}$ is a Markov process on the state space $\{0, 1, \ldots, N\}$, where $N \in \mathbb{N}_0 \cup \{\infty\}$, characterized by (nonnegative) birth rates $\{b_j, j \in \{0, 1, \ldots, N\}\}$ and (nonnegative) death rates $\{d_j, j \in \{1, \ldots, N\}\}$ and has a generator

$$(1.1) \qquad (\mathcal{A}h)(j) = b_j(h(j+1) - h(j)) - d_j(h(j) - h(j-1))$$

$$\text{with } j \in \{0, 1, \ldots, N\}.$$

It is well known that for any given probability distribution $\mu$ on $\mathbb{N}_0$ one can construct a birth–death process which has this distribution as its stationary distribution. For $N = \infty$, recurrence of the birth–death process $X(\cdot)$ is equivalent to $\sum_{n \geq 1} \frac{d_1 \cdots d_n}{b_1 \cdots b_n} = \infty$. The process $X(\cdot)$ is ergodic if and only if the process is recurrent and $c := 1 + \sum_{n \geq 1} \frac{b_0 \cdots b_{n-1}}{d_1 \cdots d_n} < \infty$. For $N < \infty$, irreducibility and hence ergodicity holds if $b_j > 0$, $j = 0, 1, \ldots, N-1$, $b_N = 0$ and $d_j > 0$, $j = 1, \ldots, N$.

In either case the stationary distribution of the ergodic process is given by $\mu(0) = 1/c$ and

$$(1.2) \qquad \mu(n) = \mu(0)\frac{b_0 \cdots b_{n-1}}{d_1 \cdots d_n}.$$

For any given probability distribution $\mu$ on $\mathbb{N}_0$ these recursive formulas give the corresponding class of birth–death processes which have $\mu$ as the stationary distribution. For the choice of a *unit per capita death rate* $d_j = j$ one simply obtains that

$$(1.3) \qquad b_j = \frac{\mu(j+1)}{\mu(j)}(j+1),$$

for $j \leq N - 1$. Here and throughout, if $N = \infty$, then by $j \leq N - 1$ we mean $j = 0, 1, \ldots$. The choice of these rates corresponds to the case where the *detailed balance condition*

$$(1.4) \qquad \mu(j)b_j = \mu(j+1)d_{j+1}, \qquad j = 0, 1, \ldots, N-1,$$

holds; see, for example, [1] and [12]. We will apply these well-known facts to the discrete Gibbs measure introduced in the next subsection.

Now might be a good time to note that not all probability distributions on $\mathbb{N}_0$ are given in closed form expressions; notable exceptions occur, for example, in compound Poisson distributions.

1.2. *Gibbs measures as stationary distributions of birth–death processes.* Gibbs measures can be viewed as stationary measures of birth–death processes, as follows.

We start with a discrete Gibbs measure $\mu$; for convenience we assume that $\mu$ has support $\operatorname{supp}(\mu) = \{0, \ldots, N\}$, where $N \in \mathbb{N}_0 \cup \{\infty\}$, so that $\mu$ is given by

$$(1.5) \qquad \mu(k) = \frac{1}{\mathcal{Z}} e^{V(k)} \frac{\omega^k}{k!}, \qquad k = 0, 1, \ldots, N,$$

for some function $V : \mathbb{N}_0 \to \mathbb{R}$. Here $\mathcal{Z} = \sum_{k=0}^{N} \exp(V(k)) \frac{\omega^k}{k!}$, and $\omega > 0$ is fixed. We assume that $\mathcal{Z}$ exists; $\mathcal{Z}$ is known as the partition function in models of statistical mechanics. We set $V(k) = -\infty$ for $k > N$. In terms of statistical mechanics, $\mu$ is a grand-canonical ensemble, $\omega$ is the activity and $V$ is the potential energy; see [16], Chapter 1.2.

The class of discrete Gibbs measures in (1.5) is equivalent to the class of all discrete probability distributions on $\mathbb{N}_0$ by the following simple identification: For a given probability distribution $(\mu(k))_{k \in \mathbb{N}_0}$ we have

$$(1.6) \quad V(k) = \log \mu(k) + \log k! + \log \mathcal{Z} - k \log \omega, \qquad k = 0, 1, \ldots, N,$$

with $V(0) = \log \mu(0) + \log \mathcal{Z}$. Hence $\mathcal{Z} = \frac{e^{V(0)}}{\mu(0)}$. The latter formula gives the possibility of proving convergence for a sequence of partition functions $\mathcal{Z}_n$ by using the convergence of the corresponding sequence $\mu_n(0)$.

However, the representation of a probability measure as a Gibbs measure is not unique. For example, the Poisson distribution with parameter $\lambda$, $Po(\lambda)$, can be written in the form (1.5) with $\omega = \lambda$, $V(k) = -\lambda$, $k \geq 0$, $\mathcal{Z} = 1$. Alternatively, we could have chosen $V(k) = 0$, $\omega = \lambda$, $\mathcal{Z} = e^{\lambda}$.

From (1.3), if we choose a unit per capita death rate $d_k = k$, and if we choose the birth rate

$$(1.7) \qquad b_k = \omega e^{V(k+1) - V(k)} = (k+1) \frac{\mu(k+1)}{\mu(k)},$$

for $k, k+1 \in \operatorname{supp}(\mu)$, then

$$(1.8) \quad (\mathcal{A}h)(k) = (h(k+1) - h(k)) \omega e^{V(k+1) - V(k)} + k(h(k-1) - h(k)),$$

for $k, k+1, k-1 \in \operatorname{supp}(\mu)$, $k \in \mathbb{N}$ [set $h(-1) = 0$], is the generator of the time-reversible birth–death process with invariant measure $\mu$.

Note that this choice of $d_k$ and $b_k$ ensures that the detailed balance condition (1.4) is satisfied. Hence we have chosen a birth–death process with generator of the type (1.1) which is easily seen to be an ergodic process. Namely, if $N = \infty$, then the recurrence of the corresponding process is given since

$$\sum_{n \geq 1} \frac{d_1 \cdots d_n}{b_1 \cdots b_n} = \sum_{n \geq 1} \frac{\mu(1)}{(n+1)\mu(n+1)} \geq \sum_{n \geq 1} \frac{\mu(1)}{n+1} = \infty.$$

For $N < \infty$ we have that $V(k) > -\infty$ for $k \in \mathrm{supp}(\mu)$, so that $b_k > 0$ for $k = 0, \ldots, N-1$, and as $V(N+1) = -\infty$, we obtain $b_N = 0$. Due to the unit per capita death rate the ergodicity follows. From $0 \in \mathrm{supp}(\mu)$ we obtain $c = 1/\mu(0) < \infty$. Hence $\mu$ is indeed the unique stationary distribution of the birth–death process.

In the development of Stein's method unit per capita is a common and useful choice for the death rate; see [2]. It is worth noting that there are modifications of this choice in [5] and [11].

To compare with the approach in Barbour [2] we reformulate the generator: let $h(k+1) - h(k) =: g(k+1)$; then (1.8) yields

$$(1.9) \quad (\mathcal{A}g)(k) = g(k+1)\omega e^{V(k+1)-V(k)} - kg(k), \qquad k = 0, 1, \ldots, N.$$

The generalization to case of arbitrary death rates $d_k$ is straightforward; we omit it here to streamline the paper.

## 2. Stein identity for Gibbs measures and bounds.
In view of the generator approach to Stein's method, for a test function $f : \mathrm{supp}(\mu) \to \mathbb{R}$ the appropriate *Stein equation* for $\mu$ given in (1.5) is

$$(2.1) \qquad\qquad (\mathcal{A}g)(j) = f(j) - \mu(f)$$

for $j \in \{0, \ldots, N\}$ and $\mathcal{A}$ the generator given by (1.8). Here,

$$\mu(f) := \sum_{k=0}^{N} f(k)\mu(k)$$

is the expectation of $f$ under $\mu$. We are interested in indicator functions $f(j) = I_{[j \in A]}$ for some $A \subset \mathrm{supp}(\mu)$. Thus if $W$ is a random variable on $\mathrm{supp}(\mu)$, we obtain

$$(2.2) \qquad\qquad E\mathcal{A}g(W) = P(W \in A) - \mu(A).$$

The right-hand side of (2.2) links in nicely with the total variation distance. Recall that for $P$ and $Q$ being probability distributions on $\mathbb{N}_0$, we define the *total variation distance* (*metric*) by

$$d_{\mathrm{TV}}(P,Q) := \sup_{A \subset \mathbb{N}_0} |P(A) - Q(A)|$$
$$= \sup_{f \in \mathcal{B}_1} |P(f) - Q(f)|$$
$$= \tfrac{1}{2} \sum_{k \in \mathbb{N}_0} |P(\{k\}) - Q(\{k\})|,$$

where $\mathcal{B}_1$ denotes the set of measurable functions $f$ with $0 \le f \le 1$. Hence bounding the left-hand side of (2.2) uniformly in $A \subset \mathbb{N}_0$ gives a bound on the total variation distance.

In the following we give a Stein characterization for $\mu$ given in (1.5) and a solution of the corresponding Stein equation.

Let $Z$ be a random variable distributed according to the Gibbs measure $\mu$ defined in (1.5) and for a function $g : \mathbb{N}_0 \to \mathbb{R}$ assume $E|Zg(Z)| < \infty$. From (1.9) we obtain the Stein characterization for $\mu$ that if $Z$ is distributed according to $\mu$, given by (1.5), then for every function $g : \mathbb{N}_0 \to \mathbb{R}$ with $E|Zg(Z)| < \infty$

$$(2.3) \qquad E\{\omega e^{V(Z+1)-V(Z)} g(Z+1) - Zg(Z)\} = 0.$$

If $f : \mathrm{supp}(\mu) \to \mathbb{R}$ is an arbitrary function, and $\mu$ is given by (1.5), then there exists a solution $g_{f,V} : \mathbb{N}_0 \to \mathbb{R}$ for (2.1) with operator as in (1.9),

$$(2.4) \quad g_{f,V}(k+1)\omega e^{V(k+1)-V(k)} - kg_{f,V}(k) = f(k) - \mu(f), \qquad k \le N;$$

see [5]. This solution $g_{f,V}$ is such that $g_{f,V}(0) = 0$, and for $j = 0, \dots, N-1$ the solution $g_{f,V}$ can be represented by recursion as

$$(2.5) \qquad g_{f,V}(j+1) = \frac{j!}{\omega^{j+1}} e^{-V(j+1)} \sum_{k=0}^{j} e^{V(k)} \frac{\omega^k}{k!} (f(k) - \mu(f))$$

$$(2.6) \qquad = -\frac{j!}{\omega^{j+1}} e^{-V(j+1)} \sum_{k=j+1}^{N} e^{V(k)} \frac{\omega^k}{k!} (f(k) - \mu(f)).$$

We may set $g_{f,V}(N+1) = 0$.

Having a suitable Stein equation for Gibbs measures and its solution at our disposal, the next step in Stein's method is to bound the increments of the solutions; it will turn out advantageous to bound the solutions themselves as well. For any function $g : \mathbb{N}_0 \to \mathbb{R}$ we define

$$\Delta g(j) := g(j+1) - g(j).$$

In applications often only bounds on the increments are needed, hence we start with these. Uniform bounds on the increments are also called *Stein* *factors* or *magic factors*. Nonuniform Stein factor bounds may yield better overall bounds on distributional distances and are therefore of particular interest. Lemma 2.1 gives such a nonuniform bound. The proof is given in [5], Lemma 2.4 and Theorem 2.1. We introduce the class of functions

$$(2.7) \qquad \mathcal{B} := \{f : \operatorname{supp}(\mu) \to [0, 1]\}$$

and we define

$$(2.8) \qquad F(k) := \sum_{i=0}^{k} \mu(i), \qquad \bar{F}(k) := \sum_{i=k}^{N} \mu(i).$$

LEMMA 2.1 (Nonuniform bounds for increments). *Assume that the death rates are unit per capita and assume that the birth rates in (1.7) fulfill, for each $k = 1, 2, \ldots, N-1$,*

$$(2.9) \qquad k \frac{F(k)}{F(k-1)} \geq \omega e^{V(k+1)-V(k)} \geq k \frac{\bar{F}(k+1)}{\bar{F}(k)}.$$

*Let $f \in \mathcal{B}$ and let $g_{f,V}$ be its solution to the Stein equation (2.4). Then, for every $j \in \{0, \ldots, N\}$,*

$$(2.10) \qquad \sup_{f \in \mathcal{B}} |\Delta g_{f,V}(j)| = \frac{1}{\omega} e^{V(j)-V(j+1)} \bar{F}(j+1) + \frac{1}{j} F(j-1).$$

*Moreover, for every $j \in \{0, \ldots, N\}$,*

$$\sup_{f \in \mathcal{B}} |\Delta g_{f,V}(j)| \leq \frac{1}{j} \wedge \frac{e^{V(j)}}{\omega e^{V(j+1)}}.$$

REMARK 2.2. Condition (2.9) is Condition (C2) in [5], Lemma 2.4. In this paper, three more conditions are formulated, which are all equivalent to (C2). For example, if the death rates are unit per capita and the birth rates are nonincreasing:

$$(2.11) \qquad e^{V(k+1)-V(k)} \leq e^{V(k)-V(k-1)}, \qquad k = 0, 1, \ldots, N,$$

Condition (C4) in [5] is satisfied and (2.10) holds.

REMARK 2.3. Reference [11] gives an elegant recursive proof of (2.10) for a choice of birth rates and death rates which make the method of exchangeable pairs work. In particular, unit per capita death rates are not used in her results.

Under a slightly weaker condition than (2.9) it is possible to derive nonuniform bounds on the solution $g_{f,V}$ of the Stein equation itself, as follows.

LEMMA 2.4.   *Consider the solution $g_{f,V}$ of the Stein equation (2.4), where $f \in \mathcal{B}$, given in (2.7). Assume that, for each $k = 1, 2, \ldots, N-1$,*

$$(2.12) \qquad \omega e^{V(k+1)-V(k)} \geq k \frac{\bar{F}(k+1)}{\bar{F}(k)}$$

*is satisfied. Then we obtain for every $j \in \{1, \ldots, N\}$ that*

$$|g_{f,V}(j)| \leq \left( \min\left\{ \ln(j), \sum_k k\mu(k) \right\} + \frac{1}{\omega} e^{V(0)-V(1)} \right) \frac{1}{\bar{F}(j)}.$$

The proof is related to ideas of [4], pages 7–8, used in Poisson approximation.

PROOF OF LEMMA 2.4.   For $j \in \{0, \ldots, N\}$ let $U_j = \{0, 1, \ldots, j\}$. We use the notation $f\mathbf{1}_A(x) = f(x)\mathbf{1}(x \in A)$. It is easy to see from (2.5) that for $f \in \mathcal{B}$ and $A \subset \{0, \ldots, N\}$, and for $j \in \{0, \ldots, N-1\}$,

$$g_{f,V}(j+1) = \frac{j!}{\omega^{j+1}} e^{-V(j+1)} \{ \mu(f\mathbf{1}_{U_j})\mu(U_j^c) - \mu(f\mathbf{1}_{U_j^c})\mu(U_j) \}$$

and hence

$$(2.13) \qquad |g_{f,V}(j+1)| \leq \frac{j!}{\omega^{j+1}} e^{-V(j+1)} \mu(U_j)\mu(U_j^c).$$

From (2.12) we have that

$$\begin{aligned}
\frac{j!}{\omega^{j+1}} e^{-V(j+1)} \mu(U_j) &= \frac{1}{\omega} \sum_{k=0}^{j} e^{V(k)-V(j+1)} \frac{\omega^{k-j} j!}{k!} \\
&\leq \sum_{k=1}^{j} \frac{\bar{F}(0)}{\bar{F}(j+1)} \frac{1}{k} + \frac{1}{\omega} e^{V(0)-V(1)} \frac{\bar{F}(0)}{\bar{F}(j+1)} \\
&\leq \left( \ln(j) + \frac{1}{\omega} e^{V(0)-V(1)} \right) \frac{1}{\bar{F}(j+1)}.
\end{aligned}$$

Alternatively,

$$\begin{aligned}
\frac{j!}{\omega^{j+1}} e^{-V(j+1)} \mu(U_j) &\leq \sum_{k=1}^{j} \frac{\bar{F}(k)}{\bar{F}(j+1)} \frac{1}{k} + \frac{1}{\omega} e^{V(0)-V(1)} \frac{\bar{F}(0)}{\bar{F}(j+1)} \\
&\leq \left( \sum_{k=1}^{N} \bar{F}(k) + \frac{1}{\omega} e^{V(0)-V(1)} \right) \frac{1}{\bar{F}(j+1)} \\
&= \left( \sum_{k=1}^{N} k\mu(k) + \frac{1}{\omega} e^{V(0)-V(1)} \right) \frac{1}{\bar{F}(j+1)},
\end{aligned}$$

proving the assertion. $\square$

We now prove a crude but usable bound for the supremum norm $\|g_{f,V}\|$ which does not require (2.9) or (2.12) to be satisfied. To this purpose we introduce the quantities

$$(2.14) \qquad \begin{aligned} \lambda_1 &:= \omega \inf_{0 \leq k \leq N-1} e^{V(k+1)-V(k)} = \inf_{0 \leq k \leq N-1} b_k, \\ \lambda_2 &:= \omega \sup_{0 \leq k \leq N-1} e^{V(k+1)-V(k)} = \sup_{0 \leq k \leq N-1} b_k. \end{aligned}$$

Note that $\frac{\lambda_2}{\lambda_1} \geq 1$ by construction.

LEMMA 2.5. *Consider the solution $g_{f,V}$ of the Stein equation (2.4), where $f \in \mathcal{B}$, given in (2.7). Assume that $\lambda_1 > 0$ and $\lambda_2 < \infty$. Then*

$$\|g_{f,V}\| \leq 2 + \frac{1}{2}\left(\frac{\lambda_2}{\lambda_1 + 1}\right)^{\lambda_2 - \lambda_1 - 2} \mathbf{1}(\lambda_2 - 2 \geq \lambda_1).$$

PROOF. As the proof follows the ideas of [4], pages 7–8, used in Poisson approximation, we only sketch it here. With the notation as for (2.13) we obtain the bounds

$$\begin{aligned} |g_{f,V}(j+1)| &\leq \frac{j!}{\omega^{j+1}} e^{-V(j+1)} \mu(U_j) \\ &= \frac{1}{\omega} \sum_{k=0}^{j} e^{V(k)-V(j+1)} \frac{\omega^{k-j} j!}{k!} \\ &\leq \lambda_1^{-1} \sum_{\ell=0}^{j} \lambda_1^{-\ell} \frac{j!}{(j-\ell)!}. \end{aligned}$$

Similarly we have

$$\begin{aligned} |g_{f,V}(j+1)| &\leq \frac{j!}{\omega^{j+1}} e^{-V(j+1)} \mu(U_j^c) \\ &= \frac{1}{\omega} \sum_{k=j+1}^{N} e^{V(k)-V(j+1)} \frac{\omega^{k-j} j!}{k!} \\ &\leq \sum_{k=j+1}^{N} \lambda_2^{k-j-1} \frac{j!}{k!} \\ &= j! e^{\lambda_2} \lambda_2^{-j-1} Po(\lambda_2)(U_j^c). \end{aligned}$$

This puts us in the situation of (1.20) and (1.21) in Chapter 1 of [4]. We obtain that for $j < \lambda_1$

$$|g_{f,V}(j+1)| \leq 2\min(1, \lambda_1^{-1/2}) \leq 2,$$

and for $j > \lambda_2 - 2$

$$|g_{f,V}(j+1)| \leq \frac{j+2}{(j+1)(j+2-\lambda_2)} \leq \frac{5}{4} < 2.$$

Thus we have proved the assertion for $j < \lambda_1$ and for $j > \lambda_2 - 2$. If $\lambda_1 > \lambda_2 - 2$, then our bound covers the whole domain.

Now assume that $S = \{\lfloor \lambda_1 \rfloor + 1, \ldots, \lfloor \lambda_2 \rfloor - 2\}$ is nonempty. For $j \in S$ we have

$$\begin{aligned}
|g_{f,V}(j+1)| &\leq \frac{j!}{\omega^{j+1}} e^{-V(j+1)} (\mu(U_j^c \cap S) + \mu(U_j^c \setminus S)) \mu(U_j) \\
&\leq 2 + \frac{j!}{\omega^{j+1}} e^{-V(j+1)} \mu(U_j^c \cap S) \\
&\leq 2 + j! \sum_{k=j+1}^{\lfloor \lambda_2 \rfloor - 2} \lambda_2^{k-j-1} \frac{1}{k!} \\
&\leq 2 + j! \lambda_2^{-j-1} e^{\lambda_2} Po(\lambda_2)\{0, \ldots, \lfloor \lambda_2 \rfloor - 2\}.
\end{aligned}$$

From [4], Proposition A.2.3(iii), page 259, the Poisson probabilities can be bounded as

$$Po(\lambda_2)\{0, \ldots, \lfloor \lambda_2 \rfloor - 2\} \leq \frac{\lambda_2}{\lambda_2 + 2 - \lfloor \lambda_2 \rfloor} Po(\lambda_2)(\lfloor \lambda_2 \rfloor - 2)$$

and so, as $\lfloor \lambda_2 \rfloor - 1 - j \geq 1$ for $j \leq \lfloor \lambda_2 \rfloor - 2$,

$$\begin{aligned}
|g_{f,V}(j+1)| &\leq 2 + \frac{1}{2} \frac{j!}{(\lfloor \lambda_2 \rfloor - 2)!} \lambda_2^{\lfloor \lambda_2 \rfloor - 2 - j} \\
&\leq 2 + \frac{1}{2}(j+1)^{-(\lfloor \lambda_2 \rfloor - 2 - j)} \lambda_2^{\lfloor \lambda_2 \rfloor - 2 - j} \\
&\leq 2 + \frac{1}{2}\left(\frac{\lambda_2}{\lambda_1 + 1}\right)^{\lfloor \lambda_2 \rfloor - 2 - j} \\
&\leq 2 + \frac{1}{2}\left(\frac{\lambda_2}{\lambda_1 + 1}\right)^{\lambda_2 - \lambda_1 - 2}.
\end{aligned}$$

This finishes the proof. $\quad\square$

REMARK 2.6. As $\lambda_1$ and $\lambda_2$ stay invariant under the reparametrization $\omega \to \tilde{\omega} = \alpha\omega$, we argue that these are reasonable quantities to employ.

REMARK 2.7.  While the examples below will show that the bound in Lemma 2.5 may not be informative, in particular examples better bounds may be obtainable in a straightforward manner. Lemma 2.5 is nevertheless useful as it gives conditions on the birth rates so that $\|g_{f,V}\|$ is bounded, and these conditions do not involve monotonicity of the birth rates.

REMARK 2.8.  Note that neither the last bound in Lemma 2.1 nor the bound in Lemma 2.5 use the normalizing constant $\mathcal{Z}$ explicitly.

Again, the generalization of the bounds to the case of arbitrary death rates $d_k$ would be straightforward.

A complication arises when we compare two distributions with nonidentical supports. Therefore it will turn out to be useful to consider the following extension from finite to infinite support for a generator $\mathcal{A}$. For convenience assume that the corresponding measure $\mu$ has $\mathrm{supp}(\mu) = \{0, 1, \ldots, n\}$, for some finite $n$, so that $\mathcal{A}$ is only defined for functions with support $\{0, 1, \ldots, n\}$ [recall that we set $g(n+1) = 0$]. We extend $\mathcal{A}$ to be defined for functions on $\{0, 1, 2, \ldots\}$ as follows:

$$
\begin{aligned}
(2.15) \qquad (\tilde{\mathcal{A}}g)(k) &:= g(k+1)\omega e^{V(k+1)-V(k)} - kg(k), \qquad k = 0, 1, \ldots, n, \\
(\tilde{\mathcal{A}}g)(k) &:= -kg(k), \qquad k \geq n+1.
\end{aligned}
$$

Thus when there are more than $n$ individuals the process is pure death. Now if $X \sim \mu$, then still $E\tilde{\mathcal{A}}f(X) = 0$ for all functions $f \in \mathcal{D}(\tilde{\mathcal{A}})$, the domain of $\tilde{\mathcal{A}}$, since the operator $\tilde{\mathcal{A}}$ represents a birth–death process with the same invariant distribution as $\mathcal{A}$; see (1.2). Next we extend the solution of the Stein equation, so that for $f : \{0, 1, \ldots, n\} \to \mathbb{R}$, the solution $g_{f,V}$ is defined by (2.5) for $k \in \{0, 1, \ldots, n\}$, and by

$$
(2.16) \qquad g_{f,V}(k) := \frac{\mu(f)}{k}, \qquad k \geq n+1.
$$

[Note that our formula (2.5) would have yielded $g_{f,V}(n+1) = 0$.] For a related suggestion see [4], Chapter 9.2. The above definition ensures that the Stein equation (2.4) is still satisfied. However, the bounds on the solution of (2.4) change slightly. In contrast to (2.7), let

$$
(2.17) \qquad \mathcal{B}_0 := \{f : \mathbb{N}_0 \to [0, 1], f(x) = 0 \text{ for } x \notin \mathrm{supp}(\mu)\}.
$$

LEMMA 2.9.  *Let $\tilde{\mathcal{A}}$, associated with $\mu, V$, and $\omega$ be given as in (2.15), where $\mu$ has $\mathrm{supp}(\mu) = \{0, 1, \ldots, n\}$. Let $f \in \mathcal{B}_0$ and let $g_{f,V}$ be the solution of the Stein equation (2.4) given in (2.5) for $k \leq n$, and as in (2.16) for* $k \geq n+1$. *Defining the sum $\sum_{l=j+1}^{n}$ as $0$ for $j \geq n$, with $\lambda_1$ and $\lambda_2$ being given in (2.14) for $\mu$,*

$$|g_{f,V}(j)| \leq 2 + \frac{1}{2}\left(\frac{\lambda_2}{\lambda_1+1}\right)^{\lambda_2-\lambda_1-2}\mathbf{1}(\lambda_2 - 2 \geq \lambda_1), \qquad j = 0, 1, \ldots, n;$$

$$|g_{f,V}(j)| \leq \frac{1}{j}, \qquad j \geq n+1.$$

*In the case that (2.12) is satisfied we also have*

$$|g_{f,V}(j)| \leq \left(\min\left\{\ln(j), \sum_k k\mu(k)\right\} + \frac{1}{\omega}e^{V(0)-V(1)}\right)\frac{1}{\bar{F}(j)}, \qquad j = 0, 1, \ldots, n,$$

*and if (2.9) is satisfied, then*

$$\sup_{f \in \mathcal{B}_l}|\Delta g_{f,V}(j)| \leq \begin{cases} \dfrac{1}{j} \wedge \dfrac{e^{V(j)}}{\omega e^{V(j+1)}}, & j \leq n, \\[2mm] \dfrac{1}{j}, & j > n. \end{cases}$$

PROOF. The first assertion follows directly from Lemma 2.5 and (2.16). The second assertion follows from Lemma 2.1 for $j \leq n-1$, and for $j \geq n+1$ it follows from (2.16) that

$$g_{f,V}(j+1) - g_{f,V}(j) = \mu(f)\left\{\frac{1}{j+1} - \frac{1}{j}\right\},$$

so that

$$|g_{f,V}(j+1) - g_{f,V}(j)| < \frac{1}{j}.$$

For $j = n$ we obtain

$$g_{f,V}(n+1) - g_{f,V}(n) = \mu(f)\frac{1}{n+1} - \left\{-\frac{1}{n}(f(n) - \mu(f))\right\}$$

$$= \frac{1}{n}\left(f(n) - \mu(f)\frac{1}{n+1}\right),$$

so that

$$|g_{f,V}(n+1) - g_{f,V}(n)| < \frac{1}{n}. \qquad \square$$

We conclude this section with some examples.

EXAMPLE 2.10 (Poisson distribution $Po(\lambda)$ with parameter $\lambda > 0$). We use $\omega = \lambda$, $V(k) = -\lambda$, $\mathcal{Z} = 1$. The Stein operator is

$$(\mathcal{A}g)(k) = g(k+1)\lambda - kg(k).$$

Applying Lemma 2.1 ($b_k = \lambda$ is nonincreasing in $k$) gives the nonuniform bound

$$(2.18) \qquad \sup_{f \in \mathcal{B}} |\Delta g_{f,\mathrm{pos}}(k)| \leq \frac{1}{k} \wedge \frac{1}{\lambda},$$

which leads to the well-known uniform bound $1 \wedge 1/\lambda$; see [4], Lemma 1.1.1. Yet [4], Lemma 1.1.1 gives the better bound $\lambda^{-1}(1 - e^{-\lambda}) \leq \min(1, 1/\lambda)$. In the Poisson case the right-hand side of (2.10) gives

$$1/k \sum_{l=0}^{k-1} e^{-\lambda} \frac{\lambda^l}{l!} + 1/\lambda \sum_{l=k+1}^{\infty} e^{-\lambda} \frac{\lambda^l}{l!}.$$

This sum can be bounded by

$$e^{-\lambda} \left( \frac{1}{\lambda} \sum_{l=1}^{\infty} \frac{\lambda^l}{l!} \right) = \frac{e^{-\lambda}}{\lambda} (e^\lambda - 1) = \frac{1}{\lambda}(1 - e^{-\lambda}),$$

and therefore we obtain for any $k \geq 1$

$$(2.19) \qquad \sup_{f \in \mathcal{B}} |\Delta g_{f,\mathrm{pos}}(k)| \leq \frac{1}{\lambda}(1 - e^{-\lambda})$$

as in [4], Lemma 1.1.1. The bounds (2.18) and (2.19) lead to

$$(2.20) \qquad \sup_{f \in \mathcal{B}} |\Delta g_{f,\mathrm{pos}}(k)| \leq \frac{1}{k} \wedge \frac{1}{\lambda}(1 - e^{-\lambda}),$$

which is a slight improvement of (2.18) and of [5].

For the Poisson distribution, $\lambda_1 = \lambda_2 = \lambda$ in Lemma 2.5, and the bound 2 from Lemma 2.5 is poor compared to the bound $\|g\| \leq \min(1, \lambda^{-1/2})$ from Lemma 1.1.1 on page 7 in [4]. The bound in Lemma 2.4 is cumbersome to compute; but for $\lambda \geq \sqrt{2}$ and for $j = 1$, say, the nonuniform bound $|g(1)| \leq (\lambda(1 - e^{-\lambda}))^{-1}$ is slightly more informative than the uniform bound $\|g\| \leq \lambda^{-1/2}$. The nonuniform bound improves with increasing $\lambda$, and deteriorates with increasing $j$.

EXAMPLE 2.11 (Binomial distribution with parameters $n$ and $0 < p < 1$). We use $\omega = \frac{p}{1-p}$, $V(k) = -\log((n-k)!)$, and $\mathcal{Z} = (n!(1-p)^n)^{-1}$. The Stein operator is

$$(\mathcal{A}g)(k) = g(k+1) \frac{p(n-k)}{(1-p)} - kg(k).$$

In [4] and [7] we find $(\mathcal{A}g)(k) = g(k+1)p(n-k) - (1-p)kg(k)$ as the operator which is equivalent to our formulation. The birth rates $(b_k)_k$ are nonincreasing and we obtain from Lemma 2.1 that

$$\sup_{f \in \mathcal{B}} |\Delta g_{f,\mathrm{bin}}(k)| \leq \frac{1}{(1-p)k} \wedge \frac{1}{p(n-k)}.$$

The proof of [4], Lemma 9.2.1, implicitly contains the same nonuniform bound. Formula (18) in [7] gives a bound of the form

$$\sup_{f \in \mathcal{B}} |\Delta g_{f,\text{bin}}(k)| \leq (1 - p^{n+1} - (1-p)^{n+1})/((n+1)(1-p)p)$$

for every $0 < k < n$. Simple calculations show that the bound (18) in [7] is for some cases better and for some cases worse in comparison to our nonuniform result. From Lemma 2.5 with $\lambda_1 = 1$, $\lambda_2 = n$ we obtain a bound on $\|g_{f,\text{bin}}\|$ which is of order $O(e^{n \ln n})$ and therefore in most applications not useful. The nonuniform bound Lemma 2.4

$$|g(j)| \leq \left(1 - \sum_{k=0}^{j-1} \binom{n}{k} p^k (1-p)^{n-k}\right)^{-1} \left(\min\{\ln(j), np\} + \frac{1-p}{np}\right)$$

will still be informative for $j$ small and $np$ large.

EXAMPLE 2.12 (Geometric distribution).   Consider $\mu(k) = p(1-p)^k$ for $k = 0, 1, \ldots$. The Stein operator is

$$(\mathcal{A}g)(k) = g(k+1)(1-p)(k+1) - kg(k).$$

The birth rates $b_k = (1-p)(k+1)$ fulfill $b_k - b_{k-1} \leq k - (k-1) = d_k - d_{k-1}$. This is condition (C4) in [5], which is sufficient for (2.9). Hence applying [5], Theorem 2.10, one obtains

$$\sup_{f \in \mathcal{B}} |\Delta g_{f,\text{geo}}(k)| \leq \frac{1}{k} \wedge \frac{1}{(1-p)(k+1)},$$

which leads to the uniform bound $1 \wedge \frac{1}{(1-p)}$.

We can improve their bounds calculating the right-hand side of (2.10) explicitly:

$$|\Delta g_{f,V}(j)| = \frac{p}{(1-p)(j+1)} \left(\frac{(1-p)^{j+1}}{p}\right) + \frac{1}{j}(1 - (1-p)^j)$$

$$= \frac{(1-p)^j}{j+1} + \frac{1}{j}(1 - (1-p)^j) = \frac{j + 1 - (1-p)^j}{j(j+1)}.$$

Obviously $|\Delta g_{f,V}(j)| \leq \frac{1}{j}$. Using Bernoulli's inequality $(1-x)^n \geq 1 - nx$ for $x < 1$ and $n \in \mathbb{N}$ it follows that

$$|\Delta g_{f,V}(j)| \leq \frac{j + pj}{j(j+1)} = \frac{1+p}{j+1}.$$

We obtain

$$\sup_{f \in \mathcal{B}} |\Delta g_{f,\text{geo}}(k)| \leq \frac{1}{k} \wedge \frac{1+p}{(k+1)},$$

which leads to the uniform bound $1 \wedge (1 + p)$.

In Lemma 2.5 we obtain $\lambda_1 = 1 - p$, $\lambda_2 = \infty$, and hence the lemma does not give informative bounds. Lemma 2.4 gives, with $\bar{F}(j) = (1 - p)^j$, the bound $|g_{f,\mathrm{geo}}(j)| \le \frac{1}{(1-p)^j}(\min\{\ln(j), \frac{1}{p}\} + \frac{1}{1-p})$. However, with (2.6) and $f \in \mathcal{B}$ we may bound directly

$$|g_{f,\mathrm{geo}}(j+1)| \le \frac{j!}{(1-p)^{j+1}(j+1)!} \sum_{k=j+1}^{\infty} (1-p)^k$$

$$= \frac{1}{j+1} \sum_{k=0}^{\infty} (1-p)^k = \frac{1}{p(j+1)},$$

and using that $g_{f,\mathrm{geo}}(0) = 0$ we obtain the improved bound $\|g_{f,\mathrm{geo}}\| \le \frac{1}{p}$. To our knowledge this bound is new.

In [13] the author considered another Stein operator. Hence the results cannot be compared. In [14] the authors considered the same Stein operator. They obtain $\sup_{f \in \mathcal{B}} |\Delta g_{f,\mathrm{geo}}(k)| \le \frac{1}{k}$. Hence [5] already improved this result.

## 3. Size-bias couplings and Gibbs measures.

3.1. *Size-bias couplings and Stein characterizations.*  Stein's method is most powerful in presence of dependence. To treat global dependence, couplings have proved a useful tool. For discrete, nonnegative random variables, so-called *size-bias* couplings fit nicely into our framework as they link in with unit per capita death rate generators.

For any random variable $X \ge 0$ with $EX > 0$ we say that a random variable $X^*$, defined on the same probability space as $X$, has the $X$-*size-biased* distribution if

$$(3.1) \qquad EXf(X) = EX E f(X^*)$$

for all functions $f$ such that both sides of (3.1) exist. Size-bias couplings $(X, X^*)$ have been studied in connection with Stein's method for Poisson approximation (see [4], e.g.) and for normal approximations (see [9], e.g.); a general framework is given in [8].

If $X$ is discrete, then for all $x$ we have $P(X^* = x) = \frac{x}{EX} P(X = x)$. This illustrates that size biasing corresponds to sampling proportional to size; the larger a subpopulation, the more likely it is to be contained in the sample.

EXAMPLE 3.1 (Poisson distribution).  If $X \sim Po(\lambda)$, then from the Stein operator in Example 2.10 we read off that $X^* = X + 1$ has the $Po(\lambda)$-size-biased distribution. In [4] a related coupling is used, namely the *reduced* size-biased coupling $(X, X_*)$, with $X_* = X^* - 1$.

EXAMPLE 3.2 (Bernoulli distribution). If $X \sim Be(p)$, it is easy to see that $X^* = 1$ has the $Be(p)$-size-biased distribution, for all $p \in (0, 1]$. As an aside, this shows that $X^*$ does not uniquely define $X$.

EXAMPLE 3.3 (Sum of nonnegative random variables). Let $X_1, \ldots, X_n$ be nonnegative, $EX_i = \mu_i > 0$, $W = \sum_{i=1}^n X_i$, $EW = \mu$, $\text{Var}(W) = 1$. Goldstein and Rinott [9] give the following construction, valid also in the presence of dependence. Choose an index $I$ from $\{1, \ldots, n\}$ according to $P(I = i) = \frac{\mu_i}{\mu}$. If $I = i$, replace $X_i$ by a variate $X_i^*$ having the $X_i$-size-biased distribution. If $X_i^* = x$, construct $\hat{X}_j$, $j \neq i$, such that

$$\mathcal{L}(\hat{X}_j, j \neq i | X_i^* = x) = \mathcal{L}(X_j, j \neq i | X_i = x).$$

Then

$$(3.2) \qquad W^* = \sum_{j \neq I} \hat{X}_j + X_I^*$$

has the $W$-size-biased distribution; see [9], construction after Lemma 2.1.

As in the Poisson case, the size-bias coupling can be used to derive a characterization of a Gibbs measure, as follows.

LEMMA 3.4. *Let $X \geq 0$ be such that $0 < E(X) < \infty$, and let $\mu$ be a discrete Gibbs measure given in (1.5). If $X \sim \mu$ and $X^*$ having the $X$-size-biased distribution is defined on the same probability space as $X$, then for every function $g : \mathbb{N}_0 \to \mathbb{R}$ such that $E|Xg(X)|$ exists,*

$$(3.3) \qquad \omega E e^{V(X+1)-V(X)} g(X+1) = \omega E e^{V(X+1)-V(X)} E g(X^*).$$

PROOF. In view of (2.3), for any $X \geq 0$ with $EX > 0$ and $(X, X^*)$ a size-biased coupling we have, for any bounded function $g : \mathbb{N}_0 \to \mathbb{R}$,

$$E\{\omega e^{V(X+1)-V(X)} g(X+1) - Xg(X)\}$$
$$= E\{\omega e^{V(X+1)-V(X)} g(X+1) - EX E g(X^*)\}.$$

For $X$ having distribution (1.5), the expectation is given by

$$EX = \omega E e^{V(X+1)-V(X)}.$$

From $EX > 0$ it follows that $\omega \neq 0$. The result now follows from the Stein characterization (2.3) of (1.5) for $\omega \neq 0$.  □

Indeed in (3.3) the factor $\omega$ cancels. However, in a moment we shall relate (3.3) to the Stein equation. In order not to get confused about the solutions for the Stein equation with and without the $\omega$ involved, we have decided to keep the $\omega$.

Lemma 3.4 provides a new formulation of the Stein approach (2.2): Let $W \geq 0$ with $0 < EW < \infty$. If $W$ has distribution $\mu$, then for all $f \in \mathcal{B}$

$$
\begin{aligned}
(3.4) \quad & Ef(W) - \mu(f) \\
& = \omega\{Ee^{V(W+1)-V(W)}g_{f,V}(W+1) - Ee^{V(W+1)-V(W)}Eg_{f,V}(W^*)\},
\end{aligned}
$$

where $g_{f,V}$, given in (2.5), is the solution of the Stein equation (2.4).

3.2. *Comparison of distributions via their generators.* From (2.2), for any random variable $W$ and any measurable set $A$ we can find a function $f = f_\mu$ such that $P(W \in A) - \mu(A) = E\mathcal{A}f(W)$, where $\mathcal{A}$ is the generator associated with the target distribution $\mu$. Applications of Stein's method usually continue by bounding the right-hand side $E\mathcal{A}f(W)$. Nevertheless, Stein's method can also be employed to compare two distributions by comparing their generators. While generators may in general not be available, for any discrete distribution we have implicitly constructed a generator in (1.9). It is from this angle that we shall bound distances between Gibbs distributions. Assume that $\mu_1$ and $\mu_2$ are two distributions with generators $\mathcal{A}_1$ and $\mathcal{A}_2$, respectively. Then, for $W \sim \mu_2$, $E\mathcal{A}_2f_{\mu_1}(W) = 0$, and therefore

$$
\mu_1(A) - P(W \in A) = E\mathcal{A}_1f_{\mu_1}(W) - E\mathcal{A}_2f_{\mu_1}(W) = E(\mathcal{A}_1 - \mathcal{A}_2)f_{\mu_1}(W).
$$

Hence we can compare two discrete Gibbs distributions by comparing their birth rates and their death rates, as follows.

THEOREM 3.5. *Let $\mu_1$ have generator $\mathcal{A}_1$ as in (1.9) and corresponding $(\omega_1, V_1)$, and let $\mu_2$ have generator $\mathcal{A}_2$, and corresponding $(\omega_2, V_2)$, both described in terms of unit per capita death rates. Suppose that $\mathcal{D}(\mathcal{A}_1) = \mathcal{D}(\mathcal{A}_2)$. Then, for $X_2 \sim \mu_2$, $f \in \mathcal{B}$, if $g_{f,V_1}$ is the solution of the Stein equation for $\mu_1$,*

$$
\begin{aligned}
& \left| Ef(X_2) - \int f \, d\mu_1 \right| \\
& \quad \leq \min\Bigg\{ \|g_{f,V_1}\| E(X_2) \\
(3.5) \quad & \qquad \times \left( \frac{|\omega_1 - \omega_2|}{\omega_2} + \frac{\omega_1}{\omega_2} E|e^{(V_1(X_2^*)-V_1(X_2^*-1))-(V_2(X_2^*)-V_2(X_2^*-1))} - 1| \right), \\
& \qquad \|g_{f,V_2}\| E(X_1) \\
& \qquad \times \left( \frac{|\omega_2 - \omega_1|}{\omega_1} + \frac{\omega_2}{\omega_1} E|e^{(V_2(X_1^*)-V_2(X_1^*-1))-(V_1(X_1^*)-V_1(X_1^*-1))} - 1| \right) \Bigg\}.
\end{aligned}
$$

PROOF. For $X_2 \sim \mu_2$, $f \in \mathcal{B}$, if $g_{f,V_1}$ is the solution of the Stein equation for $\mu_1$, we have $Ef(X_2) - \mu_1(f) = E\mathcal{A}_1 g_{f,V_1}(X_2) = E(\mathcal{A}_1 - \mathcal{A}_2)g_{f,V_1}(X_2)$. Using the size biasing and (3.3) we obtain

$$
\begin{aligned}
&Ef(X_2) - \mu_1(f) \\
&\quad = Eg_{f,V_1}(X_2+1)(\omega_1 e^{V_1(X_2+1)-V_1(X_2)} - \omega_2 e^{V_2(X_2+1)-V_2(X_2)}) \\
&\quad = \omega_1 Eg_{f,V_1}(X_2+1)e^{V_2(X_2+1)-V_2(X_2)}e^{V_1(X_2+1)-V_1(X_2)-(V_2(X_2+1)-V_2(X_2))} \\
&\qquad - E(X_2)Eg_{f,V_1}(X_2^*) \\
&\quad = \frac{\omega_1}{\omega_2}E(X_2)Eg_{f,V_1}(X_2^*)e^{(V_1(X_2^*)-V_1(X_2^*-1))-(V_2(X_2^*)-V_2(X_2^*-1))} \\
&\qquad - E(X_2)Eg_{f,V_1}(X_2^*) \\
&\quad = \frac{\omega_1 - \omega_2}{\omega_2}E(X_2)Eg_{f,V_1}(X_2^*) \\
&\qquad + \frac{\omega_1}{\omega_2}E(X_2)Eg_{f,V_1}(X_2^*)\{e^{(V_1(X_2^*)-V_1(X_2^*-1))-(V_2(X_2^*)-V_2(X_2^*-1))} - 1\}.
\end{aligned}
$$

Taking absolute values together with the triangle inequality, and observing that the argument is symmetric in the indices 1 and 2, we obtain (3.5). $\quad\square$

REMARK 3.6. Theorem 3.5 has the advantage that the partition function $\mathcal{Z}$ may not be needed to assess the distance.

REMARK 3.7. Note that it is not surprising that $\|g_{f,V_1}\|$ appears in the bound; even if we compare two Poisson distributions, with different means $\lambda_1$ and $\lambda_2$, the bound would be of the type $|\lambda_1 - \lambda_2|\|g_{f,V_1}\|$.

When $\mathcal{D}(\mathcal{A}_1) \neq \mathcal{D}(\mathcal{A}_2)$, we first extend $\mathcal{A}_1$ to $\tilde{\mathcal{A}}_1$, defined for functions on $\{0, 1, 2, \dots\}$ as in (2.15), and apply then Lemma 2.9. For convenience assume that $\operatorname{supp}(\mu_1) = \{0, 1, \dots, n\}$ and $\operatorname{supp}(\mu_2) = \{0, 1, 2, \dots\}$. Using similar calculations as for Theorem 3.5, we derive the following result.

COROLLARY 3.8. Let $\mu_1$ have generator $\mathcal{A}_1$ and corresponding $(\omega_1, V_1)$, and let $\mu_2$ have generator $\mathcal{A}_2$, and corresponding $(\omega_2, V_2)$, both described in terms of unit per capita death rates. Suppose that $\operatorname{supp}(\mu_1) = \{0, 1, \dots, n\}$ and that $\operatorname{supp}(\mu_2) = \{0, 1, 2, \dots\}$. Then, for $X_2 \sim \mu_2$ and $f \in \mathcal{B}_0$ as in (2.17), if $g_{f,V_1}$ is the solution of the Stein equation for $\mu_1$ as in (2.15),

$$
\begin{aligned}
&\left| Ef(X_2) - \int f \, d\mu_1 \right| \\
&\qquad \leq \min\Big\{ \|g_{f,V_1}\|E(X_2)
\end{aligned}
$$

$$\times \left( \frac{|\omega_1 - \omega_2|}{\omega_2} + \frac{\omega_1}{\omega_2} E \big| e^{(V_1(X_2^*) - V_1(X_2^*-1)) - (V_2(X_2^*) - V_2(X_2^*-1))} - 1 \big| \right),$$

(3.6)
$$\|g_{f,V_2}\| E(X_1)$$

$$\times \left( \frac{|\omega_2 - \omega_1|}{\omega_1} + \frac{\omega_2}{\omega_1} E \big| e^{(V_2(X_1^*) - V_2(X_1^*-1)) - (V_1(X_1^*) - V_1(X_1^*-1))} - 1 \big| \right) \Big\}$$

$$+ \mu_1(f) \sum_{k=n+1}^{\infty} \mu_2(k).$$

The extra term $\mu_1(f) \sum_{k=n+1}^{\infty} \mu_2(k)$ in the bound (3.6) as compared to (3.5) accounts for the extension $\tilde{\mathcal{A}}_1$ of the generator $\mathcal{A}_1$.

PROOF OF COROLLARY 3.8. We proceed as in the proof of Theorem 3.5. Let $\tilde{\mathcal{A}}_1$ be the extension of $\mathcal{A}_1$ as in (2.15). For $X_2 \sim \mu_2$, $f \in \mathcal{B}_l$, if $g_{f,V_1}$ as in (2.16) is the solution of the Stein equation for $\mu_1$, we have

$$Ef(X_2) - \mu_1(f)$$
$$= E\tilde{\mathcal{A}}_1 g_{f,V_1}(X_2)$$
$$= E(\tilde{\mathcal{A}}_1 - \mathcal{A}_2) g_{f,V_1}(X_2)$$
$$= E g_{f,V_1}(X_2+1)(\omega_1 e^{V_1(X_2+1) - V_1(X_2)} - \omega_2 e^{V_2(X_2+1) - V_2(X_2)})$$
$$= \sum_{k=0}^{n-1} P(X_2=k) g_{f,V_1}(k+1)(\omega_1 e^{V_1(k+1) - V_1(k)} - \omega_2 e^{V_2(k+1) - V_2(k)})$$
$$\quad - \sum_{k=n}^{\infty} P(X_2=k) \frac{\mu_1(f)}{k+1} \omega_2 e^{V_2(k+1) - V_2(k)},$$

where we used (2.16). Now the summand

$$\sum_{k=0}^{n-1} P(X_2=k) g_{f,V_1}(k+1)(\omega_1 e^{V_1(k+1) - V_1(k)} - \omega_2 e^{V_2(k+1) - V_2(k)})$$

can be treated exactly as in Theorem 3.5; from Lemma 2.9 we have the same bounds on the solution of the Stein equation when $k \leq n$. For the second summand we have

$$\sum_{k=n}^{\infty} P(X_2=k) \frac{\mu_1(f)}{k+1} \omega_2 e^{V_2(k+1) - V_2(k)}$$

$$= \omega_2 \mu_1(f) \sum_{k=n}^{\infty} \frac{1}{\mathcal{Z}} e^{V_2(k)} \frac{\omega_2^k}{k!(k+1)} e^{V_2(k+1) - V_2(k)}$$

$$= \mu_1(f) \sum_{k=n+1}^{\infty} \mu_2(k).$$

Thus we obtain (3.6).   $\square$

Note that the bound in Corollary 3.8 contains the measure $\mu_2$ explicitly; hence examples will require the calculation of the normalizing constant $\mathcal{Z}_2$. We now treat the instructive example that the two generators $\mathcal{A}_1$ and $\mathcal{A}_2$ have the same birth and death rates, but live on different domains.

EXAMPLE 3.9.   Suppose that the distributions $\mu_1$ and $\mu_2$ have supports $\mathrm{supp}(\mu_1) = \{0, 1, \ldots, n\}$ and $\mathrm{supp}(\mu_2) = \{0, 1, 2, \ldots\}$, respectively. Let $\mu_1$ have generator $\mathcal{A}_1$ and corresponding $(\omega, V)$, and let $\mu_2$ have generator $\mathcal{A}_2$, and corresponding same $(\omega, V)$, both described in terms of unit per capita death rates, so that the two generators have the same birth rates and death rates for $k = 0, \ldots, n$. Then, for $X_2 \sim \mu_2$ and $f \in \mathcal{B}_0$, if $g_{f,V_1}$ is the solution of the Stein equation for $\mu_1$, it follows from (3.6) that

$$\left| Ef(X_2) - \int f \, d\mu_1 \right| \leq \mu_1(f) \sum_{k=n+1}^{\infty} \mu_2(k).$$

In particular it follows that

$$d_{\mathrm{TV}}(\mu_1, \mu_2) \leq \sum_{k=n+1}^{\infty} \mu_2(k),$$

stating that the difference between the two distributions can be bounded by the total mass of the domain, under the respective distribution, which is in the support of the one distribution, but not of the other distribution.

Indeed the bound is sharp, as from (1.2) we have that $\mu_1(k) = \alpha \mu_2(k)$ for $k = 0, 1, \ldots, n$, with $\alpha = (\sum_{k=0}^{n} \mu_2(k))^{-1} > 1$. Hence

$$d_{\mathrm{TV}}(\mu_1, \mu_2) = \tfrac{1}{2} \sum_{k=0}^{n} (\alpha - 1)\mu_2(k) + \tfrac{1}{2} \sum_{k=n+1}^{\infty} \mu_2(k) = \sum_{k=n+1}^{\infty} \mu_2(k).$$

Theorem 3.5 and Corollary 3.8 invite the study of the total variation distance of discrete probability distributions with different supports of the corresponding operators, for example considering the distance between a binomial and a hypergeometric distribution. As this leads away from the flow of the paper, we do not pursue it here. Instead we turn to the application of our approach to Gibbs measures arising in interacting particle systems.

**4. Application to lattice approximation in statistical physics.** Now we apply our results to the model studied by Chayes and Klein [6]. Their setup is as follows. Assume $A \subset \mathbb{R}^d$ is a rectangle; denote its volume by $|A|$. Consider the intersection of $A$ with the $d$-dimensional lattice $n^{-1}\mathbb{Z}^d$. For each site $m$ in this intersection we associate a Bernoulli random variable $X_m^n$ which takes value 1, with probability $p_m^n$, if a particle is present at $m$ and 0 otherwise. If the collection $(X_m^n)_m$ is independent, the joint distribution can be interpreted as the Gibbs distribution for an ideal gas on the lattice. The Poisson convergence theorem states that, for $n$ going to infinity, when preserving the average density of particles the lattice ideal gas distribution converges weakly to a Poisson distribution, which is the standard Gibbs distribution for an ideal gas in the continuum. On physical grounds one expects that a similar result might hold for interacting particles. The model is as follows. Pick $n \in \mathbb{N}$ and suppose that $A$ can be partitioned into a regular array of $d(n)$ sub-rectangles $\{S_1^n, \ldots, S_{d(n)}^n\}$ with volumes $v(S_m^n) = \frac{z_m^n}{z_n}$. Here, $z_m^n > 0$ for each $m, n$. Thus

$$|A| = \frac{1}{z_n} \sum_{m=1}^{d(n)} z_m^n.$$

As a guideline motivated by [6] we choose $z_m^n > 0$ and $z_n$ such that $z_m^n \to 0$ and $z_n \to z > 0$ for $n \to \infty$. For each $m, n$ choose a point $q_m^n \in S_m^n$. Now the following class of functions is considered.

ASSUMPTION 4.1. *Let $(f_k)_k$ be a sequence of functions satisfying $f_0 \equiv 1$ and for each $k \geq 1$, $f_k(x_1, \ldots, x_k)$ is a nonnegative function, Riemann integrable on $A^k$, such that $f_k$ is a symmetric function for each $k$; that is, for any permutation $\sigma$, $f_k(x_{\sigma(1)}, \ldots, x_{\sigma(k)}) = f_k(x_1, \ldots, x_k)$.*

Let $X_m^n, 1 \leq m \leq d(n)$, be 0–1 random variables with joint density function given as follows. If $a_1, \ldots, a_{d(n)}$ are such that $a_i \in \{0, 1\}$, and if

$$k = \sum_{m=1}^{d(n)} a_m = \sum_{m=1}^{k} a_{i_m}$$

is the sum of the nonzero $a_i$'s, then

$$P(X_1^n = a_1, \ldots, X_{d(n)}^n = a_{d(n)})$$

$$= K(d(n))^{-1} f_k(q_{i_1}^n, \ldots, q_{i_k}^n) \prod_{m=1}^{d(n)} (z_m^n)^{a_m},$$

where the normalizing constant is

$$K(d(n)) = \sum_{\mathbf{a} \in \{0,1\}^{d(n)}} \sum_k \mathbf{1}_{\{\sum a_i = k\}} f_k(q_{i_1}^n, \ldots, q_{i_k}^n) \prod_{m=1}^{d(n)} (z_m^n)^{a_m}.$$

Define
$$S_n = \sum_{m=1}^{d(n)} X_m^n.$$

Then, due to the symmetry of $f_k$,

$$(4.1) \quad P(S_n = k) = \frac{\frac{(z_n)^k}{k!} \sum_{i_1=1}^{d(n)} \cdots \sum_{i_k=1}^{d(n)} f_k(q_{i_1}^n, \ldots, q_{i_k}^n) \prod_{m=1}^{k} v(S_{i_m}^n)}{\sum_{k=0}^{d(n)} \frac{(z_n)^k}{k!} \sum_{i_1=1}^{d(n)} \cdots \sum_{i_k=1}^{d(n)} f_k(q_{i_1}^n, \ldots, q_{i_k}^n) \prod_{m=1}^{k} v(S_{i_m}^n)}.$$

Note that we can write (4.1) as a Gibbs measure $\mu_n$ of the form (1.5) with

$$(4.2) \quad \omega_n = z_n, \qquad V_n(k) = \log(W_n(k)), \qquad k \in \{0, 1, \ldots \, d(n)\},$$

where

$$W_n(k) = \sum_{i_1=1}^{d(n)} \cdots \sum_{i_k=1}^{d(n)} f_k(q_{i_1}^n, \ldots, q_{i_k}^n) \prod_{m=1}^{k} v(S_{i_m}^n), \qquad k \in \{0, \ldots, d(n)\}.$$

Let $S$ be a nonnegative integer-valued random variable defined by

$$(4.3) \quad P(S = k) = \frac{\frac{z^k}{k!} \int_{A^k} f_k(x_1, \ldots, x_k) \, dx_1 \cdots dx_k}{\sum_{k=0}^{\infty} \frac{z^k}{k!} \int_{A^k} f_k(x_1, \ldots, x_k) \, dx_1 \cdots dx_k}.$$

We write (4.3) as a Gibbs measure $\mu$ of the form (1.5) with

$$\omega = z, \qquad V(k) = \log(W(k)), \qquad k \in \{0, 1, \ldots\},$$

where

$$W(k) = \int_{A^k} f_k(x_1, \ldots, x_k) \, dx_1 \cdots dx_k, \qquad k \in \{0, 1, \ldots\}.$$

In [6] the following class of functions is considered. Let $(f_k)_k$ be a sequence of functions satisfying Assumption 4.1 and, in addition:

(a) There exists a constant $C$ such that $f_k(x_1, \ldots, x_k) \leq C^k$ for all $k \geq 0$.

(b) $f_k(x_1, \ldots, x_k) = 0$ if $x_i = x_j$ for some $i \neq j$. This is not a necessary condition; see [6].

In [6] it is shown that under these additional conditions, $S_n \Rightarrow S$, that is, $S_n$ converges weakly to $S$, under the conditions that $\lim_{n\to\infty} z_n = z$ and that

$$\lim_{n\to\infty} \left( \max_{1 \leq m \leq d(n)} z_m^n \right) = 0.$$

For the proof of this convergence, condition (a) on $(f_k)_k$ enters to ensure that the Riemann sum converges, so that the normalizing constants $\mathcal{Z}_n$ and $\mathcal{Z}$ for $\mu_n$ and for $\mu$, respectively, are finite. Condition (b) avoids some measure-theoretic considerations. In [6] a bound on the rate of convergence is given for the special case that $f_k = 1$, all $k \geq 0$, resulting in a Poisson approximation.

Here we give a general bound on the distance, not only in the case $f_k = 1$, $k \geq 0$.

Note that in [6] the generalization of Poisson convergence allows the authors to develop a lattice-to-continuum theory of classical statistical mechanics; similarly our results could be applied to obtain bounds on such lattice-to-continuum theory.

Using Corollary 3.8 with $\mu_1 = \mu_n$ and $\mu_2 = \mu$, we obtain

$$\left| E f(S_n) - \int f \, d\mu \right|$$
$$\leq \min \left\{ E(S_n) \|g_{f,V}\| \left( \frac{|z_n - z|}{z_n} + \frac{z}{z_n} E \left| \frac{W(S_n^*) W_n(S_n^* - 1)}{W(S_n^* - 1) W_n(S_n^*)} - 1 \right| \right), \right.$$
$$\left. E(S) \|g_{f,V_n}\| \left( \frac{|z - z_n|}{z} + \frac{z_n}{z} E \left| \frac{W_n(S^*) W(S^* - 1)}{W_n(S^* - 1) W(S^*)} - 1 \right| \right) \right\}$$
$$+ \mu_n(f) \sum_{k=d(n)+1}^{\infty} \mu(k).$$

Writing out the size-biased distribution gives the following result.

THEOREM 4.2.  *Let $S_n, \mu, z_n$ be as above. Let us assume that $(f_k)_k$ satisfies Assumption 4.1 and condition (*a*). Then*

$$d_{\mathrm{TV}}(\mathcal{L}(S_n), \mu)$$
$$\leq \max(\|g_{f,V}\|, \|g_{f,V_n}\|)$$
$$\times \min \left\{ \left( E(S_n) \frac{|z_n - z|}{z_n} + \frac{z}{z_n} \sum_k k \mu_n(k) \left| \frac{W(k) W_n(k-1)}{W(k-1) W_n(k)} - 1 \right| \right), \right.$$
$$\left. \left( E(S) \frac{|z - z_n|}{z} + \frac{z_n}{z} \sum_k k \mu(k) \left| \frac{W_n(k) W(k-1)}{W_n(k-1) W(k)} - 1 \right| \right) \right\}$$
$$+ \mu_n(f) \sum_{k=d(n)+1}^{\infty} \mu(k).$$

REMARK 4.3.  Note that for Theorem 4.2, no monotonicity assumptions on the birth and death rates in the corresponding birth–death process are needed. For the bound to be informative, though, we need that $\max(\|g_{f,V}\|, \|g_{f,V_n}\|) < \infty$; conditions to ensure this behavior are given in Lemma 2.5.

For the case of nonincreasing birth rates and unit per capita death rates, we also obtain bounds based on the increments $\Delta g$. Recall the notation of Example 3.3.

PROPOSITION 4.4.   *Let $S_n$ be the sum of $X_1, \ldots, X_n$, where $X_i \sim Be(p_i)$, and construct $S_n^*$ via (3.2), and let $\mu$ be given by (1.5). Assume that the birth rates in (1.7) are nonincreasing, and assume that the death rates are unit per capita. Then*

$$d_{\mathrm{TV}}(\mathcal{L}(S_n), \mu)$$

$$\leq \omega E\{e^{V(S_n+1)-V(S_n)}$$

(4.4)
$$\times \min\{S(S_n, S_n^*-1),$$

$$(|S_n - S_n^* + 1|/\omega)e^{V(\min(S_n, S_n^*-1))-V(\min(S_n, S_n^*-1)+1)}\}\}$$

$$+ \|g_{f,V}\|\omega E|e^{V(S_n+1)-V(S_n)} - E e^{V(S_n+1)-V(S_n)}|,$$

*where we put $S(x, y) = \sum_{\ell=\min(x,y)+1}^{\max(x,y)} \frac{1}{\ell}$.*

PROOF.   With the notation as in Example 3.3, denote the conditional expectation given $I = i$ by $E_i$. Let $f \in \mathcal{B}$. If $X_i \sim Be(p_i)$, $i = 1, \ldots, n$, then $X_i^* = 1$. Put

$$\hat{S}_{n,i} = \sum_{j \neq i} \hat{X}_j.$$

By (3.4) and (3.2), we have

$$|Ef(S_n) - \mu(f)|$$

$$= \omega \left| \sum_{i=1}^{n} \frac{p_i}{\sum_{j=1}^{n} p_j} \left\{ E_i[e^{V(S_n+1)-V(S_n)} g_{f,V}(S_n+1)] \right. \right.$$

$$\left. \left. - E e^{V(S_n+1)-V(S_n)} E g_{f,V}\left(\sum_{j \neq i} \hat{X}_j + 1\right) \right\} \right|$$

$$= \omega \left| \sum_{i=1}^{n} \frac{p_i}{\sum_{j=1}^{n} p_j} E_i[e^{V(S_n+1)-V(S_n)})(g_{f,V}(S_n+1) - g_{f,V}(\hat{S}_{n,i}+1))] \right.$$

$$\left. + E_i g_{f,V}(\hat{S}_{n,i}+1)(e^{V(S_n+1)-V(S_n)} - E e^{V(S_n+1)-V(S_n)}) \right|$$

$$\leq \frac{\omega}{\sum_{j=1}^{n} p_j} \sum_{i=1}^{n} p_i E_i\{e^{V(S_n+1)-V(S_n)}$$

$$\times \min\{S(S_n, \hat{S}_{n,i}),$$

$$(|S_n - \hat{S}_{n,i}|/\omega)$$

$$\times e^{V(\min(S_n, \hat{S}_{n,i})) - V(\min(S_n, \hat{S}_{n,i}) + 1)}\} | S_n \neq \hat{S}_{n,i}\}$$

$$\times P(S_n \neq \hat{S}_{n,i})$$

$$+ \|g_{f,V}\| \omega E |e^{V(S_n + 1) - V(S_n)} - E e^{V(S_n + 1) - V(S_n)}|,$$

where we used Lemma 2.1 for the last inequality.   □

Note that in contrast to Theorem 4.2, we did not employ the difference between generators in the proof.

4.1. *Poisson approximation for a sum of indicators.* In the Poisson case we can use the bound (2.20) on $\Delta g$ instead of Lemma 2.1 to obtain improved overall bounds, as follows. Let $X_i \sim Be(p_i)$, $i = 1, \ldots, n$; then as in Example 3.3 we have for $S_n = \sum_{j=1}^{n} X_j$ that $S_n^* = \hat{S}_{n,I} + 1$, where as above $\hat{S}_{n,i} = \sum_{j \neq i} \hat{X}_j$. Let $\lambda = ES_n = \sum_{j=1}^{n} p_j$. Recall that if $\mu = Po(\lambda)$, then $\omega = \lambda$ and $V(k) = -\lambda$. Hence using (4.5) we obtain

$$|Ef(S_n) - Po(\lambda)(f)| = \left| \sum_{i=1}^{n} p_i E_i(g_{f,V}(S_n + 1) - g_{f,V}(\hat{S}_{n,i})) \right|.$$

Calculating as in the proof of Proposition 4.4 and using (2.20) we get

$$d_{\mathrm{TV}}(\mathcal{L}(S_n), Po(\lambda))$$

$$(4.5) \qquad \leq \sum_{i=1}^{n} p_i E_i \left\{ \min\left\{ S(S_n, \hat{S}_{n,i}), \frac{|S_n - \hat{S}_{n,i}|}{\lambda}(1 - e^{-\lambda}) \right\} \Big| S_n \neq \hat{S}_{n,i} \right\}$$

$$\times P(S_n \neq \hat{S}_{n,i}).$$

In particular we recover the bound for couplings in [4], Theorem 1.B,

$$(4.6) \qquad d_{\mathrm{TV}}(\mathcal{L}(S_n), Po(\lambda)) \leq \frac{1 - e^{-\lambda}}{\lambda} \sum_{i=1}^{n} p_i E_i |S_n - \hat{S}_{n,i}|.$$

EXAMPLE 4.5. Assume that the variables $S_n$ and $S_n^*$ can be coupled such that

$$(4.7) \qquad |S_n - \hat{S}_{n,i}| \leq 1, \qquad i = 1, \ldots, n.$$

Then we can improve the bound (4.6), using

$$E_i\{S(S_n, \hat{S}_{n,i}) | S_n \neq \hat{S}_{n,i}\} P(S_n \neq \hat{S}_{n,i})$$

$$\leq E_i \left\{ \frac{S_n - \hat{S}_{n,i}}{1 + \hat{S}_{n,i}} \Big| S_n \neq \hat{S}_{n,i} \right\} P(S_n \neq \hat{S}_{n,i})$$

$$\leq P(S_n \neq \hat{S}_{n,i}) \left[ P\{\hat{S}_{n,i} = 0 | S_n \neq \hat{S}_{n,i}\} + \frac{1}{2} P\{\hat{S}_{n,i} \geq 1 | S_n \neq \hat{S}_{n,i}\} \right]$$

$$= P(S_n \neq \hat{S}_{n,i}) \frac{1}{2} (1 + P\{\hat{S}_{n,i} = 0 | S_n \neq \hat{S}_{n,i}\}).$$

EXAMPLE 4.6.   In case that the $X_i$'s are independent and we have $\hat{S}_{n,i} = S_n - X_i$, so that (4.7) is satisfied, and (4.6) yields that

$$(4.8) \qquad d_{\mathrm{TV}}(\mathcal{L}(S_n), Po(\lambda)) \leq \sum_{i=1}^{n} p_i^2 \frac{1 - e^{-\lambda}}{\lambda}.$$

The bound (4.8) coincides with the bound (1.23) given in [4]. We can improve on it by using that

$$P\{\hat{S}_{n,i} = 0 | S_n \neq \hat{S}_{n,i}\} = P\left( \sum_{j \neq i} X_j = 0 | X_i = 1 \right) = \prod_{j \neq i} (1 - p_j).$$

This results in the bound for the independent case

$$d_{\mathrm{TV}}(\mathcal{L}(S_n), Po(\lambda)) \leq \sum_{i=1}^{n} p_i^2 \min \left\{ \frac{1}{2} \left( 1 + \prod_{j \neq i} (1 - p_j) \right), \frac{1 - e^{-\lambda}}{\lambda} \right\}.$$

To see when the above bound is an improvement on the bound (1.23) given in [4], consider the inequality

$$\frac{1}{2}(1 + x) < \frac{1 - e^{-\lambda}}{\lambda}.$$

We rearrange this inequality to give

$$\lambda x < 2 - \lambda - 2e^{-\lambda}.$$

Numerical calculations yield that the right-hand side is positive when $0 < \lambda < 1.59362$. For such $\lambda$ and some $y$ such that $0 < y < 2 - \lambda - 2e^{-\lambda}$ we can construct strongly inhomogeneous cases such that

$$\prod_{j \neq i} (1 - p_j) \leq \lambda^{-1} y$$

for some indices $i$. Examples are, for $\lambda = 1$, the vector $\mathbf{p} = (p_1, \ldots, p_n) = (1, 0, 0, \ldots, 0)$, or the vector $\mathbf{p} = (p_1, \ldots, p_n) = (1 - 1/(n-1), 1/(n-1)^2, 1/(n-1)^2, \ldots, 1/(n-1)^2)$ for $n \geq 3$. In such cases, our bound provides an improvement.

Example 4.6 can only be an improvement on (4.6) when the underlying random variables are not identically distributed; in the independent and *identically distributed* case, we recover exactly the bound given in [4], (1.23), page 8.

4.2. *An example with repelling interaction.* Our first example with interaction partitions $A = [0,1]$ into $n$ intervals, the $i$th interval given by $S_i^n = [\frac{i-1}{n}, \frac{i}{n}]$, and we choose $q_i^n = \frac{2i-1}{2n}$, the midpoint of the interval $S_i^n$. Then each $S_i^n$ has volume $\frac{1}{n} = \frac{\lambda}{\lambda n}$; we choose $z_i^n = \frac{\lambda}{n}$ and $z_n = \lambda$. Note the freedom of choice in $\lambda > 0$. We consider the set of functions $f_k$ given by $f_0 = 1, f_1(x) = 1$ for all $x$, and for $k \geq 2$, $f_k(x_1, \ldots, x_k) = 0$ if $x_i = x_j$ for some $i \neq j$, and otherwise

$$f_k(x_1, \ldots, x_k) = \sum_{1 \leq i \neq j \leq k} (x_i - x_j)^2,$$

so that values $x_1, \ldots, x_n$ which are far away from each other are preferred. Then clearly $f_k$ satisfies Assumption 4.1 as well as conditions (a) and (b). To avoid trivialities we assume that $n \geq 2$.

In our setup we choose 0–1 random variables $X_m^n$, $1 \leq m \leq n$, with density function

$$P(X_1^n = a_1, \ldots, X_n^n = a_n) \propto n^{-k} \lambda^k \sum_{1 \leq i \neq j \leq k} (q_i^n - q_j^n)^2 \qquad \text{if } k = \sum_{m=1}^{n} a_m \geq 2;$$

$$P(X_1^n = a_1, \ldots, X_n^n = a_n) \propto \frac{\lambda}{n} \qquad \text{if } k = 1;$$

$$P(X_1^n = 0, \ldots, X_n^n = 0) \propto 1.$$

With $S_n = \sum_{m=1}^{n} X_m^n$ we have with (4.1) that

$$P(S_n = k) \propto \frac{\lambda^k}{k!} W_n(k),$$

where

$$W_n(k) = n^{-k} \sum_{i_1=1}^{n} \cdots \sum_{i_k=1}^{n} f_k(q_{i_1}^n, \ldots, q_{i_k}^n), \qquad k = 0, \ldots, n.$$

Let $S$ be a nonnegative integer-valued random variable defined as in (4.3) by

$$P(S = k) = \frac{\frac{\lambda^k}{k!} \int_{[0,1]^k} f_k(x_1, \ldots, x_k) \, dx_1 \cdots dx_k}{\sum_{k=0}^{\infty} \frac{\lambda^k}{k!} \int_{[0,1]^k} f_k(x_1, \ldots, x_k) \, dx_1 \cdots dx_k} \propto \frac{\lambda^k}{k!} W(k),$$

with

$$W(k) = \int_{[0,1]^k} f_k(x_1, \ldots, x_k) \, dx_1 \cdots dx_k, \qquad k \geq 0.$$

It is immediate that $W(0) = W(1) = 1$, and for $k \geq 2$,

$$W(k) = \sum_{i=1}^{k} \sum_{j \neq i} \int_0^1 \int_0^1 (x_i - x_j)^2 \, dx_i \, dx_j = \frac{k(k-1)}{6}.$$

Thus our limiting Gibbs measure $\mu$ has normalizing constant

$$\mathcal{Z} = 1 + \lambda + \frac{1}{6} \sum_{k \geq 2} \frac{\lambda^k}{(k-2)!} = 1 + \lambda + \frac{\lambda^2}{6} e^{\lambda}$$

and is therefore given by

$$\mu(0) = \left(1 + \lambda + \frac{\lambda^2}{6} e^{\lambda}\right)^{-1},$$

$$\mu(1) = \frac{\lambda}{1 + \lambda + \frac{\lambda^2}{6} e^{\lambda}},$$

$$\mu(k) = \left(1 + \lambda + \frac{\lambda^2}{6} e^{\lambda}\right)^{-1} \frac{\lambda^k}{(k-2)!}, \qquad k \geq 2.$$

To assess the distance between the distributions of $S_n$ and of $S$ we note that $W_n(0) = W_n(1) = 1 = W(0) = W(1)$, and some algebra yields, for $k \geq 2$, that

$$W_n(k) = n^{-k} \sum_{i_1=1}^{n} \cdots \sum_{i_k=1}^{n} \sum_{l \neq s} \left(\frac{2i_s - 1}{2n} - \frac{2i_j - 1}{2n}\right)^2$$

$$= n^{-k-2} \sum_{1 \leq l \neq s \leq k} \sum_{i=1}^{n} \sum_{j=1}^{n} (i-j)^2 n^{k-2}$$

$$= 2n^{-4} k(k-1) \left\{ \sum_{i=1}^{n} i^2 - \left(\sum_{i=1}^{n} i\right)^2 \right\}$$

$$= \frac{k(k-1)(n+1)(n-1)}{6n^2}.$$

Thus

$$\frac{W_n(k)}{W_n(k-1)} = \frac{W(k)}{W(k-1)}$$

for $k = 0, \ldots, n$, and we are in the situation of Example 3.9. Therefore it follows that

$$d_{\mathrm{TV}}(\mu_n, \mu) = \sum_{k=n+1}^{\infty} \mu(k)$$

$$= \left(1 + \lambda + \frac{\lambda^2}{6} e^{\lambda}\right)^{-1} \sum_{k=n+1}^{\infty} \frac{\lambda^k}{k!}$$

$$\leq \frac{\lambda^{n+1}e^\lambda}{(n+1)!(1+\lambda+\frac{\lambda^2}{6}e^\lambda)}.$$

Summarizing, we have:

COROLLARY 4.7. *Let $S_n$ and $S$ be as constructed above. Then*

$$d_{\mathrm{TV}}(\mathcal{L}(S_n), \mu) \leq \frac{\lambda^{n+1}e^\lambda}{(n+1)!(1+\lambda+\frac{\lambda^2}{6}e^\lambda)}.$$

4.3. *A second example with interaction.* Using the same partition of $A = [0,1]$ and the same notation as in the previous example, and choosing for simplicity $\omega = \lambda = 1$, we now consider the set of functions $(f_k)_k$ given by $f_0(x) = f_1(x) = 1$ for all $x$, and for $k \geq 2$,

$$f_k(x_1, \ldots, x_k) = \prod_{1 \leq i \neq j \leq k} x_i x_j.$$

Clearly $(f_k)_k$ satisfy Assumption 4.1 as well as conditions (a) and (b). We take $q_i^n = \frac{i-1}{n}$, the left endpoint of the interval $S_i^n$. Assume that $n \geq 3$. Along the lines of the calculations in the previous example we obtain now $W(0) = W(1) = 1$ and for $k \geq 2$

$$W(k) = k^{-k},$$

so that the distribution of $S$ given in (4.3) is

$$\mu(k) = \frac{1}{\mathcal{Z}} \frac{1}{k!} k^{-k}, \qquad k \geq 0.$$

We note that

$$1 \leq \mathcal{Z} = 1 + \sum_{k=1}^{\infty} \frac{1}{k!} e^{-k\ln k} \leq e^{1/e}.$$

For the distribution of $S_n$ given in (4.1) we obtain $W_n(0) = W_n(1) = 1$, and for $k \geq 2$,

$$W_n(k) = n^{-k^2} \left(\sum_{i=0}^{n-1} i^{k-1}\right)^k.$$

We employ an integral approximation, obtaining

$$\left(\frac{n-1}{n}\right)^{k^2} k^{-k} < W_n(k) < k^{-k}, \qquad k \geq 2.$$

Observe that for $k = 0, \ldots, n$,

$$\left(\frac{n-1}{n}\right)^{k^2} < \frac{W_n(k)}{W(k)} < 1.$$

For Theorem 4.2, we may thus bound for $k = 2, \ldots, n$,

$$\left| \frac{W(k)W_n(k-1)}{W(k-1)W_n(k)} - 1 \right|$$

$$\leq \max\left\{ 1 - \left(\frac{n-1}{n}\right)^{k^2}; \left(\frac{n}{n-1}\right)^{k^2} - 1 \right\}$$

$$= \left(\frac{n}{n-1}\right)^{k^2} - 1.$$

Now we calculate

$$\sum_k k\mu_n(k) \left| \frac{W(k)W_n(k-1)}{W(k-1)W_n(k)} - 1 \right| \leq \sum_k k\mu_n(k) \left\{ \left(\frac{n}{n-1}\right)^{k^2} - 1 \right\}.$$

Taylor's expansion gives the inequality

$$0 \leq \left(\frac{n}{n-1}\right)^{k^2} - 1 \leq \frac{k^2}{n-1}\left(1 + \frac{1}{n-1}\right)^{k^2} < \frac{k^2}{n-1}e^{k^2/(n-1)} < \frac{k^2}{n-1}e^k,$$

for $k \leq n$ and $n \geq 3$. Applying this inequality, we obtain that

$$\sum_k k\mu_n(k) \left| \frac{W(k)W_n(k-1)}{W(k-1)W_n(k)} - 1 \right| \leq \frac{1}{(n-1)\mathcal{Z}_n} \sum_{k=1}^{n-1} \frac{1}{k!}e^{-k\log(k)+k}$$

$$< \frac{e^e}{(n-1)\mathcal{Z}_n} < \frac{2e^e}{n},$$

where we used $\mathcal{Z}_n = 1 + \sum_{k=1}^{\infty} \frac{1}{k!}W_n(k) \geq 1$ in the last inequality. Next we estimate

$$\sum_{k=n+1}^{\infty} \mu(k) = \frac{1}{\mathcal{Z}} \sum_{k=n+1}^{\infty} \frac{1}{k!}k^{-k} \leq \sum_{k=n+1}^{\infty} \frac{1}{k!}n^{-k} \leq e^{1/n}\frac{n^{-(n+1)}}{(n+1)!},$$

where we used that $\mathcal{Z} > 1$.

Note that $\lambda_2 = \sup_k \frac{(k+1)^{k+1}}{k^k} = \infty$, so that we cannot apply Lemma 2.5; instead we find an alternative bound on $\|g_{f,V}\|$ as follows. From (2.6) we have

$$|g_{f,V}(j+1)| \leq \frac{j!}{\omega^{j+1}}e^{-V(j+1)} \sum_{k=j+1}^{N} e^{V(k)}\frac{\omega^k}{k!} \leq \sum_{k=0}^{\infty} \frac{1}{k!}k^{-k} = \mathcal{Z} \leq e^{1/e}.$$

COROLLARY 4.8. *Let $S_n$ and $S$ be as constructed above, and let $n \geq 3$. Then*

$$d_{\mathrm{TV}}(\mathcal{L}(S_n), \mu) \leq 2\frac{e^{e+1/e}}{n} + e^{1/n}\frac{n^{-(n+1)}}{(n+1)!}.$$

**Acknowledgments.** The authors would like to thank Hans Zessin for pointing out the work of Preston as well as the work of Chayes and Klein. We would also like to thank Andrew Barbour and two anonymous referees for helpful comments.

Fakultät für Mathematik
Ruhr-Universität Bochum
NA 3/68
D-44780 Bochum
Germany
E-mail: peter.eichelsbacher@ruhr-uni-bochum.de

Department of Statistics
University of Oxford
1 South Parks Road
Oxford OX1 3TG
United Kingdom
E-mail: reinert@stats.ox.ac.uk